\documentclass[preprint,number,sort&compress]{elsarticle}

\journal{\ }

\usepackage{amsmath}
\usepackage{amsfonts}
\usepackage{mathtools}

\usepackage{xcolor}

\newtheorem{thm}{Theorem}
\newtheorem{lem}[thm]{Lemma}
\newtheorem{conj}{Conjecture}
\newdefinition{rmk}{Remark}
\newdefinition{exa}{Example}
\newproof{pf}{Proof}
\newcommand{\proofofref}{}
\newproof{zpot}{Proof of \proofofref}
\newenvironment{pot}[1]%
	{\renewcommand{\proofofref}{#1}\zpot}%
 {\endzpot}

\begin{document}

\begin{frontmatter}

	\title{On the separation of solutions to fractional differential equations of order $\alpha \in (1,2)$}
	
	\author[1]{Renu Chaudhary}
	\ead{renu.chaudhary@thws.de}
	
	\author[1]{Kai Diethelm\corref{cor1}}
	\ead{kai.diethelm@thws.de}
	
	\author[1]{Safoura Hashemishahraki}
	\ead{safoura.hashemishahraki@thws.de}
	
	\affiliation[1]{organization={Faculty of Applied Natural Sciences and Humanities, 
						Technical University of Applied Sciences Würzburg-Schweinfurt},
			addressline={Ignaz-Schön-Str.~11},
			postcode = {97421},
			city = {Schweinfurt},
			country = {Germany}
			}
			
	\cortext[cor1]{Corresponding author}

	\begin{abstract}
		Given the Caputo-type fractional
		differential equation $D^\alpha y(t) = f(t, y(t))$ with $\alpha \in (1, 2)$, 
		we consider two distinct solutions $y_1, y_2 \in C[0,T]$ to this equation subject to 
		different sets of initial conditions.
		In this framework, we discuss nontrivial upper and lower bounds for the difference 
		$|y_1(t) - y_2(t)|$ for $t \in [0,T]$.
		The main emphasis is on describing how such bounds are related to the differences of the
		associated initial values.
	\end{abstract}

	
\end{frontmatter}

\section{Introduction}
\label{sec:intro}

In a recent paper \cite{DT:separation1}, one of the authors contributed to the development of
a theory addressing the following question:
\begin{quote}
	\em
	Given some $\alpha \in (0,1]$ and a fractional initial value problem
	\[
		D^\alpha y(t) = f(t, y(t)), \qquad
		y(0) = y_0,
	\]
	with differential operators of Caputo type and starting point $0$ 
	and order $\alpha$,
	by how much does the solution change if the initial value $y_0$ is replaced by, say, $y_0 + \delta$
	with some $\delta \ne 0$?
\end{quote}
Upper bounds for the changes in the solution can be easily derived using a Gronwall argument \cite{DF},
but in many practically relevant cases these can be found to be excessively large \cite[p.~170]{DT:separation1}.
Therefore it was the primary goal of \cite{DT:separation1} to provide much more rigorous upper bounds and 
to complement these with correspondingly strict lower bounds. In this work here, we now extend this theory 
to the case of differential equations of order $\alpha \in (1,2)$.

Thus, we will concentrate on the initial value problems
\begin{equation}
	\label{eq:ivp}
	D^\alpha y_k(t) = f(t, y_k(t)), \qquad
	y_k(0) = y_{k,0}, \quad
	y_k'(0) = y_{k,1}
	\qquad
	(k = 1, 2)
\end{equation}
where $\alpha \in (1,2)$, i.e.\ we look at two functions $y_1$ and $y_2$ 
that solve \emph{the same fractional differential equation} but 
satisfy \emph{different initial conditions}, and derive strict upper and lower bounds for the difference
$y_1(t) - y_2(t)$ for $t \in [0,T]$ with some $T > 0$. 
Throughout our work, $D^\alpha$ will denote the Caputo differential operator of order $\alpha$ 
with starting point $0$ \cite[Definition 3.2]{Diethelmbook}.

Throughout this paper, we use the standard notation 
\[
	E_{\alpha, \beta}(z) = \sum_{k=0}^\infty \frac {z^k}{\Gamma(\alpha k + \beta)}
\]
to denote the two-parameter Mittag-Leffler function with parameters $\alpha >0$ and $\beta > 0$
\cite{GKMR2020}. With this convention, we can state our first result.

\begin{thm}
	Assume that $\alpha \in (1,2)$ and that the function 
	$f : [0,T] \times \mathbb R \to \mathbb R$ is continuous and satisfies a 
	Lipschitz condition with respect to its second argument whose Lipschitz constant is
	denoted by $L$. Then, unique solutions $y_1$ and 
	$y_2$ to the two initial value problems given in \eqref{eq:ivp} exist on $[0,T]$. Moreover,
	for all $t \in [0, T]$ we have that 
	\begin{equation}
		\label{eq:gronwall}
		|y_1(t) - y_2(t)| \le \left( |y_{1,0} - y_{2,0}| + t |y_{1,1} - y_{2,1}| \right) E_{\alpha,1}(L t^\alpha).
	\end{equation}
\end{thm}

\begin{pf}
	The existence and uniqueness of the solutions $y_1$ and $y_2$ follows from classical 
	standard results, see, e.g., \cite[Theorem 6.4]{Tisdell2012}.
	A careful inspection of the proof of the corresponding Gronwall inequality in this case 
	\cite[Theorem 6.20]{Diethelmbook} (see also \cite{DF}) then yields the inequality \eqref{eq:gronwall}.
	\qed
\end{pf}

The estimate \eqref{eq:gronwall} has two fundamental weaknesses: First, it does not provide 
a nontrivial lower bound for $|y_1(t) - y_2(t)|$ at all. And second, by definition, the Lipschitz constant 
$L$ is always positive (except for the trivial case when $f(t,y)$ depends only on $t$). This implies that, when the
interval $[0,T]$ is large, the Mittag-Leffler function needs to be evaluated for a large positive argument, 
and thus in view of its known asymptotic behaviour 
\[
	E_{\alpha, 1}(L t^\alpha) = \alpha^{-1} \exp(L^{1/\alpha} t) + O(t^{-\alpha}) 
	\qquad \text{ for } t \to \infty
\]
(see \cite[Theorem 4.3]{GKMR2020}), the bound on the right-hand side of \eqref{eq:gronwall} becomes
very large. In Section \ref{sec:separation}, we will derive refined estimates that address both these weaknesses.

For the case $\alpha \in (0,1)$, it has been shown in \cite{DU:shooting1} that results of the type to 
be developed below
are useful in designing and analyzing numerical methods for solving so-called terminal value
problems for fractional differential equations of order $\alpha$. In the same spirit, we believe that
our results described here can be exploited to construct numerical methods for fractional boundary
value problems of the form
\begin{equation}
	\label{eq:bvp}
	D^\alpha y(t) = f(t, y(t)), \qquad
	y(0) = y_0, \quad y(T) = y_1,
\end{equation}
with $\alpha \in (1,2)$ and some $T > 0$. We intend to elaborate on this topic in a forthcoming paper
\cite{Di:shooting-bvp}.

It will become evident below that, when transferring the results of \cite{DT:separation1}
from the case $\alpha \in (0,1]$ treated there to $\alpha \in (1,2)$ that we discuss here, 
some steps of the proofs can be carried over directly,
but in some other respects, significant changes are necessary. Specifically, the following issues require
major modifications:
\begin{itemize}
\item In the case $\alpha \in (0,1]$, the initial value problem needs one initial condition to assert the
	well-posedness; in our case $\alpha \in (1, 2)$ here, two such conditions are required.
	Thus the methods need to be refined in order to allow incorporating the second initial condition.
\item When looking at different initial conditions, the case $\alpha \in (0,1]$ necessarily implies that
	we have two solutions that do not coincide at the initial point. This leads to a substantial simplification
	of the proof that we do not always have if $\alpha \in (1,2)$ because in this case 
	we may have two solutions with different first derivatives but identical function values at the
	initial point.
\item An inhomogeneous linear fractional differential equation of order $\alpha \in (0,1]$ with constant coefficients
	has a solution given (via the variation-of-constants formula) in terms of the 
	Mittag-Leffler functions $E_{\alpha, 1}$ and $E_{\alpha,\alpha}$. For $\alpha \in (0,1)$, both these 
	functions are strictly positive on the entire real line, see \cite{HAV2013}. 
	This property is essential in the proofs of \cite{DT:separation1}. In the case $\alpha \in (1,2)$, 
	the solutions comprise Mittag-Leffler functions $E_{\alpha, 1}$, $E_{\alpha, 2}$ and 
	$E_{\alpha, \alpha}$, and for the values of $\alpha$ that are applicable now, some
	or all (depending on the precise value of $\alpha$) of these function do have changes of sign
	on $\mathbb R$. This fact demands the proof techniques to be adapted to a significant extent.
\end{itemize}

\section{The separation of the solutions}
\label{sec:separation}

Much as in \cite{DT:separation1}, we will break up the analysis into multiple steps. 
Specifically, we will begin by showing some auxiliary results needed to prove the main theorems. 
Next, we will discuss the linear case and then finally investigate the general problem
involving nonlinear differential equations. In each case, we will keep one of the two initial values
(function value and first derivative of the solution to the initial value problem) fixed and estimate 
the effects of changes of the other one, and afterwards interchange the roles of the two initial values. The overall 
result then follows by combining the partial results in an appropriate way.

\subsection{Auxiliary results}
\label{subs:aux}

Our first result is a comparison statement for solutions to weakly singular Volterra equations.
Unlike the other results of our paper, it is not restricted to the case $\alpha \in (1,2)$ but it
holds for general $\alpha > 0$.

\begin{lem}
	\label{lem:comparison}
	Let $J = [0,T]$ with some $T > 0$ or $J = [0, \infty)$.
	Moreover, let $\alpha > 0$ and assume that the continuous functions $v, v_1, w ,w_1 : J \to \mathbb R$ 
	and $g : J \times \mathbb R \to \mathbb R$ satisfy 
	\[
		v(t) \le v_1(t) + \frac 1 {\Gamma(\alpha)} \int_0^t (t-s)^{\alpha-1} g(s, v(s)) \,  \mathrm d s
	\]
	and
	\[
		w(t) \ge w_1(t) + \frac 1 {\Gamma(\alpha)} \int_0^t (t-s)^{\alpha-1} g(s, w(s)) \,  \mathrm d s,
	\]
	respectively, for all $t \in J$. Suppose further that $g(t, x)$ is nondecreasing in $x$ for each
	fixed $t \in J$ and that $v_1(t) < w_1(t)$ for all $t \in J$. Then $v(t) < w(t)$ for all $t \in J$.
\end{lem}

\begin{pf}
	This has been claimed and proven for the case $\alpha \in (0,1)$ by Cong and Tuan
	\cite[Lemma 3.4]{CT2017} (see also \cite[Theorem 2.1]{LV2008} for a slightly modified version). 
	A close inspection of the proof given in \cite{CT2017} reveals that all arguments used there are valid 
	for $\alpha \ge 1$ too.
	\qed
\end{pf}

Next, we generalize the so-called variation-of-constants formula developed for the case
$0 < \alpha \le 1$ in \cite[Lemma 3.1]{CT2017} to the case $1 < \alpha <2$ that is
of interest here.

\begin{thm}
	\label{thm:voc}
	Let $\alpha \in (1,2)$,  $y_{1,0}, y_{1,1}, M \in \mathbb R$, and $T > 0$. Assume that the
	function $g : [0, T] \times \mathbb R \to \mathbb R$ is continuous and satisfies a Lipschitz condition
	with respect to the second variable. Moreover, let the function $y_1 : [0,T] \to \mathbb R$ be the solution
	to the initial value problem
	\[
		D^\alpha y_1(t) = M y_1(t) + g(t, y_1(t)), 
		\qquad y_1(0) = y_{1, 0} \mbox{ and } y_1'(0) = y_{1,1}.
	\]
	Then, for all $t \in [0,T]$, we have
	\begin{align}
		\label{eq:var-of-const}
		y_1(t) & = y_{1,0} E_{\alpha, 1}(M t^\alpha) + y_{1,1} t E_{\alpha, 2}(M t^\alpha) \\
		\nonumber
			& \qquad {} 
			 + \int_0^t (t - \tau)^{\alpha-1} E_{\alpha, \alpha}(M (t-\tau)^\alpha) g(\tau, y_1(\tau)) \, \mathrm d \tau .
	\end{align}
\end{thm}

\begin{pf}
	The proof is completely analog to the proof of the result for the case $\alpha \in (0,1]$ given in
	\cite[Lemma 3.1]{CT2017}.
	\qed
\end{pf}

\begin{rmk}
	In fact, it is easy to see that the result of Theorem \ref{thm:voc} can be generalized even more
	without essential changes in the proof. Specifically, we may admit arbitrary $\alpha > 0$. Although 
	we do not need this generalization in the context of the present work, we mention it here for the
	sake of completeness:

	\emph{Let $\alpha > 0$, $m = \lceil \alpha \rceil$, $y_{1,\ell}, M \in \mathbb R$ ($\ell = 0, 1, \ldots, m-1$), 
	and $T > 0$. Assume that the
	function $g : [0, T] \times \mathbb R \to \mathbb R$ is continuous and satisfies a Lipschitz condition
	with respect to the second variable. Moreover, let the function $y_1 : [0,T] \to \mathbb R$ be the solution
	to the initial value problem
	\[
		D^\alpha y_1(t) = M y_1(t) + g(t, y_1(t)), 
		\qquad D^\ell y_1(0) = y_{1, \ell} \quad (\ell = 0, 1, \ldots, m-1).
	\]
	Then, for all $t \in [0,T]$, we have
	\[
		y_1(t)  = \sum_{\ell = 0}^{m - 1} y_{1,\ell} t^\ell E_{\alpha, \ell+1}(M t^\alpha)
			 + \int_0^t (t - \tau)^{\alpha-1} E_{\alpha, \alpha}(M (t-\tau)^\alpha) g(\tau, y_1(\tau)) \, \mathrm d \tau .
	\]
	}
\end{rmk}

Based on Lemma \ref{lem:comparison} and Theorem \ref{thm:voc}, we can now prove a first separation theorem for solutions to
the differential equation from \eqref{eq:ivp}. This result generalizes the findings of \cite[Theorem 3.5]{CT2017}
to the case $\alpha \in (1,2]$ on sufficiently small intervals. We shall see in Example \ref{ex:meet} below that this restriction is 
necessary in the sense that---in contrast to the situation encountered for $0 < \alpha \le 1$---one can 
construct counterexamples if the intervals are too large.

\begin{thm}
	\label{thm:no-intersect}
	Consider the two initial value problems stated in \eqref{eq:ivp} on some interval $[0, T]$ with $\alpha \in (1,2)$
	and the associated solutions $y_1$ and $y_2$ subject to the following assumptions:
	\begin{enumerate}[(a)]
	\item The function $f$ is continuous and 
		satisfies a Lipschitz condition with respect to the second variable with Lipschitz constant $L$.
	\item The quantity $T^* \in (0, T]$ is chosen such that 
		\begin{subequations}
			\label{eq:tstar}
		\begin{equation}
			\label{eq:tstar1}
			(y_{2,0} - y_{1,0}) E_{\alpha,1}(-L t^\alpha) + (y_{2,1} - y_{1,1}) t E_{\alpha,2}(-L t^\alpha) > 0
		\end{equation}
		and
		\begin{equation}
			\label{eq:tstar2}
			E_{\alpha, \alpha}(-L t^\alpha) > 0
		\end{equation}
		\end{subequations}
		for all $t \in [0, T^*]$.
	\end{enumerate}
	Then, we have $y_1(t) < y_2(t)$ for any $t \in [0, T^*]$.
\end{thm}

\begin{rmk}
	Note that Theorem \ref{thm:no-intersect} is a qualitative result: Under the given conditions, it asserts  
	that the graphs of two solutions to the same differential equation remain separated from each other and do not 
	touch or intersect. It does, however, not quantify the distance between the solutions. We will develop
	such quantitative results in Subsections \ref{subs:linear} and \ref{subs:nonlinear} below.
\end{rmk}

\begin{rmk}
	\label{rmk:cond-signs}
	\begin{enumerate}[(a)]
	\item In condition (b) of Theorem \ref{thm:no-intersect} we demand, among others, that \eqref{eq:tstar1}
		holds for all $t \in [0, T^*]$. By looking concretely at the point $t=0$ and taking into consideration
		that $E_{ \alpha,1}(0) = 1$, this in particular means that	we must have $y_{1,0} < y_{2,0}$.
		For certain applications, in particular the one discussed in \cite{Di:shooting-bvp}, this is a little bit 
		too restrictive. Specifically, in that case we only have $y_{2,0} = y_{1,0}$ and $y_{2,1} > y_{1,1}$
		which means that \eqref{eq:tstar1} holds only for $t \in (0,T^*]$ but not for $t = 0$.
		To handle this case, we need a slightly different result that, in turn, requires a small modification of the
		proof, see Theorem \ref{thm:no-intersect-0} below.
	\item A sufficient condition for \eqref{eq:tstar1} to hold is that 
		\begin{subequations}
			\label{eq:tstar1a}
			\begin{align}
				y_{1,0} < y_{2,0}, & \quad  y_{1,1} \le y_{2,1}, \\
				\label{eq:tstar1a2}
				E_{\alpha,1}(-L t^\alpha) & > 0 \mbox{ for all } t \in [0, T^*] \\
			\intertext{and}
				\label{eq:tstar1a3}
				E_{\alpha,2}(-L t^\alpha) & > 0 \mbox{ for all } t \in [0, T^*].
			\end{align}
		\end{subequations}
	\item Since $E_{\alpha,1}$, $E_{\alpha, \alpha}$ and $E_{\alpha,2}$ are continuous functions
		with positive function values at the origin, it is clearly possible to find a value 
		$T^* > 0$ which satisfies the assumptions \eqref{eq:tstar} of Theorem \ref{thm:no-intersect}
		if $y_{1,0} < y_{2,0}$. Similarly, we can also see that it is possible to find some $T^* > 0$
		that satisfies \eqref{eq:tstar2}, \eqref{eq:tstar1a2} and \eqref{eq:tstar1a3}.
		We shall discuss this matter in more detail in Section \ref{sec:ml-zeros}.
	\end{enumerate}
\end{rmk}

\begin{pot}{Theorem \ref{thm:no-intersect}}
	Assume that the statement is not true. Since $y_1$ and $y_2$ are continuous functions
	and $y_1(0) = y_{1,0} < y_{2,0} = y_2(0)$, there exists some $t_1 \in (0, T^*]$ such that
	$y_1(t_1) = y_2(t_1)$ and that $y_1(t) < y_2(t)$ for all $t \in [0, t_1)$. Then, define
	$g(t, y) := f(t, y) + L y$ for all $t \in [0, T]$ and $y \in \mathbb R$. By Theorem \ref{thm:voc},
	we obtain for all these values of $t$ that
	\begin{subequations}
		\label{eq:y1y2voc}
		\begin{align}
			y_1(t) &= E_{\alpha, 1}(-L t^\alpha) y_{1,0} + t E_{\alpha,2}(-L t^\alpha) y_{1,1} \\
					& \qquad \nonumber
					{} + \int_0^t (t - \tau)^{\alpha-1} E_{\alpha, \alpha}(-L (t - \tau)^\alpha) 
								g(\tau, y_1(\tau)) \, \mathrm d \tau \\
		\intertext{and}
			y_2(t) &= E_{\alpha, 1}(-L t^\alpha) y_{2,0} + t E_{\alpha,2}(-L t^\alpha) y_{2,1} \\
					& \qquad \nonumber
					{} + \int_0^t (t - \tau)^{\alpha-1} E_{\alpha, \alpha}(-L (t - \tau)^\alpha) 
								g(\tau, y_2(\tau)) \, \mathrm d \tau.
		\end{align}
	\end{subequations}
	Moreover, for all these $t$ and all 
	$z, \tilde z \in \mathbb R$ with $z > \tilde z$, we observe
	\begin{equation}
		\label{eq:g-monotone}
		g(t, z) - g(t, \tilde z) 
		= f(t, z) - f(t, \tilde z) + L (z - \tilde z)
		\ge -L (z - \tilde z) + L (z - \tilde z) 
		= 0
	\end{equation}
	in view of the Lipschitz property of $f$ and therefore, by \eqref{eq:tstar2},
	\[
		 E_{\alpha, \alpha}(-L (t - \tau)^\alpha) (g(t, z) - g(t, \tilde z)) \ge 0
	\]
	for these values of $t$, $z$ and $\tilde z$ and for $\tau \in [0,t]$. 
	In other words, the function $E_{\alpha, \alpha}(-L (t - \tau)^\alpha) g(t, \cdot)$ is nondecreasing
	for any fixed $t \in [0,T]$. This observation allows us to invoke Lemma \ref{lem:comparison}
	which, in combination with \eqref{eq:tstar1} and \eqref{eq:y1y2voc}, tells us that
	$y_1(t) < y_2(t)$ for all $t \in [0, t_1]$ which contradicts our basic observation that $y_1(t_1) = y_2(t_1)$.
	\qed
\end{pot}

The following example demonstrates that the limitation on the magnitude of the quantity $T^*$ 
given in eq.~\eqref{eq:tstar} is indeed necessary.

\begin{figure}[htp]
	\centering
	\includegraphics[width=0.68\textwidth]{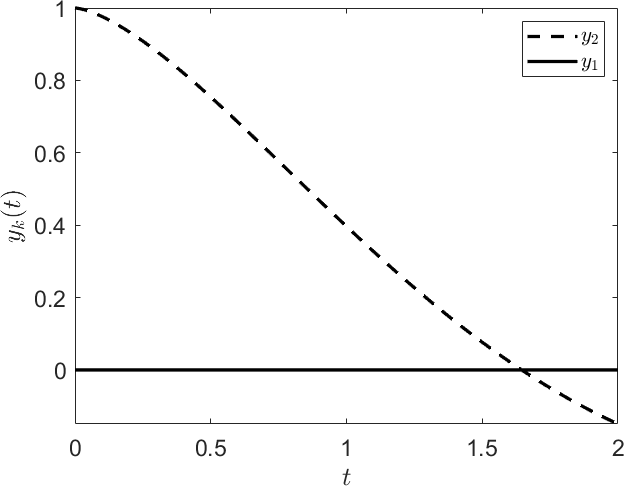}
	\caption{\label{fig:meet} The graphs of the functions $y_1(t) = 0$ and $y_2(t) = E_{1.5, 1}(-t^{1.5})$ 
		considered 	in Example \ref{ex:meet}.}
\end{figure}

\begin{exa}
	\label{ex:meet}
	Consider the differential equation $D^{1.5} y(t) = - y(t)$ which clearly satisfies the assumptions
	of Theorem \ref{thm:no-intersect}. Evidently, the solution to this equation subject to the
	initial condition $y_1(0) = y_1'(0) = 0$ is $y_1(t) = 0$. Also, it is well known 
	\cite[Theorem 7.2]{Diethelmbook} that the solution with initial conditions 
	$y_2(0) = 1$ and $y_2'(0) = 0$ is $y_2(t) = E_{1.5, 1}(-t^{1.5})$. Evidently, these
	two solutions remain separated from each other on any interval $[0, T^*]$ where $T^*$ 
	is smaller than the first positive zero of $y_2$ but not when $T^*$ is larger than this zero
	(which is located at approximately $1.645$, see Fig.\ \ref{fig:meet}), so even though both solutions 
	exist on $[0, \infty)$, their graphs remain separated from each other only on a subinterval 
	that extends from the initial point to a certain finite value.
\end{exa}

We conclude this subsection with the modification of Theorem \ref{thm:no-intersect}
that we had already announced in Remark \ref{rmk:cond-signs}(a).

\begin{thm}
	\label{thm:no-intersect-0}
	Consider the two initial value problems stated in \eqref{eq:ivp} on some interval $[0, T]$ with $\alpha \in (1,2)$
	and the associated solutions $y_1$ and $y_2$ subject to the following assumptions:
	\begin{enumerate}[(a)]
	\item The function $f$ is continuous and 
		satisfies a Lipschitz condition with respect to the second variable with Lipschitz constant $L$.
	\item $y_{1,0} = y_{2,0}$ and $y_{1,1} < y_{2,1}$.
	\item The quantity $T^* \in (0, T]$ is chosen such that the inequalities \eqref{eq:tstar2}
		and \eqref{eq:tstar1a3} hold for all $t \in [0, T^*]$.
	\end{enumerate}
	Then, we have $y_1(t) < y_2(t)$ for any $t \in (0, T^*]$.
\end{thm}

As in Remark \ref{rmk:cond-signs}(c), we also emphasize here that a value $T^*$ with the properties required in
condition (c) exists.

\begin{pf}
	By \cite[Theorem 6.25]{Diethelmbook}, we know that $y_1, y_2 \in C^1[0, T]$. Our assumption (b) 
	about the initial values asserts that $y_1'(0) < y_2'(0)$, and hence, in view
	of the continuity of $y_1'$ and $y_2'$, there exists some $\delta > 0$ such that $y_1'(t) < y_2'(t)$
	for all $t \in [0, \delta)$. Combining this with the information that $y_1(0) = y_2(0)$ (which also
	follows from our assumption about the initial values), we see that 
	\begin{equation}
		\label{eq:proof-ni0-1}
		y_1(t) < y_2(t) 
		\quad \mbox{ for all } t \in (0, \delta).
	\end{equation}
	
	Now assume that the statement of the theorem is false. Then, in view of~\eqref{eq:proof-ni0-1},
	there exists some $t_1 \ge \delta > 0$ such that 
	\begin{equation}
		\label{eq:proof-ni0-2}
		y_1(t_1) = y_2(t_1)
		\quad \mbox{ and } \quad
		y_1(t) < y_2(t)
		\mbox{ for } t \in (0, t_1).
	\end{equation}
	
	Defining the function $g$ as in the proof of Theorem \ref{thm:no-intersect}, we notice 
	again---as in~\eqref{eq:g-monotone}---that this function is strictly increasing with respect to its 
	second variable. Exploiting this monotonicity, Theorem \ref{thm:voc} and conditions (b) and (c) 
	of our Theorem, we thus derive
	\begin{eqnarray*}
		y_1(t_1) 
		& = & y_{1,0} E_{\alpha,1}(-L t_1^\alpha) + y_{1,1} t_1 E_{\alpha,2}(-L t_1^\alpha) \\
		&& {} + \int_0^{t_1} (t_1 - \tau)^{\alpha-1} E_{\alpha, \alpha}(-L (t_1 - \tau)^\alpha) 
							g(\tau, y_1(\tau)) \, \mathrm d \tau \\
		& < & y_{2,0} E_{\alpha,1}(-L t_1^\alpha) + y_{2,1} t_1 E_{\alpha,2}(-L t_1^\alpha) \\
		&& {} + \int_0^{t_1} (t_1 - \tau)^{\alpha-1} E_{\alpha, \alpha}(-L (t_1 - \tau)^\alpha) 
							g(\tau, y_1(\tau)) \, \mathrm d \tau \\
		& \le & y_{2,0} E_{\alpha,1}(-L t_1^\alpha) + y_{2,1} t_1 E_{\alpha,2}(-L t_1^\alpha) \\
		&& {} + \int_0^{t_1} (t_1 - \tau)^{\alpha-1} E_{\alpha, \alpha}(-L (t_1 - \tau)^\alpha) 
							g(\tau, y_2(\tau)) \, \mathrm d \tau \\
		& = & y_2(t_1)
	\end{eqnarray*}
	which contradicts \eqref{eq:proof-ni0-2}. Therefore, the hypothesis that our claim be false cannot hold.
	\qed
\end{pf}

\subsection{Main results for the linear case}
\label{subs:linear}

We are now in a position to derive the desired quantitative statements about the separation of
different solutions to the same fractional differential equation for the linear case, thus extending
the theory developed for $\alpha \in (0,1]$ in \cite[Theorems 4 and 5]{DT:separation1}
to the case $\alpha \in (1,2)$.

\begin{thm}
	\label{thm:sep-linear}
	Assume the hypotheses of Theorem \ref{thm:no-intersect} or Theorem \ref{thm:no-intersect-0}
	and choose $T^*$ in the corresponding way.
	Under the condition that $f(t, y) = a(t) y$
	with some continuous function $a : [0, T] \to \mathbb R$, define
	\[
		a_*(t) := \min_{\tau \in [0,t]} a(\tau)
		\quad \mbox{ and } \quad
		a^*(t) := \max_{\tau \in [0,t]} a(\tau)
	\]
	for $t \in [0,T]$. Then, 
	for all $t \in [0, T^*]$ we have
	\begin{eqnarray}
		\nonumber
		\lefteqn{|y_2(0) - y_1(0)| \cdot E_{\alpha,1}(a_*(t) t^\alpha)
				+ |y_2'(0) - y_1'(0)| \cdot t E_{\alpha,2}(a_*(t) t^\alpha)} \\
		& \le & |y_2(t) - y_1(t)| 
					\label{eq:sep-linear} \\
		& \le & |y_2(0) - y_1(0)| \cdot E_{\alpha,1}(a^*(t) t^\alpha)
				+ |y_2'(0) - y_1'(0)| \cdot t E_{\alpha,2}(a^*(t) t^\alpha) .
				\nonumber
	\end{eqnarray}
\end{thm}

\begin{pf}
	For $t \in [0, T^*]$, we define $u(t) := y_2(t) - y_1(t)$. We then see that $u$ is the unique solution 
	to the initial value problem
	\begin{equation}
		\label{eq:ivp-u}
		D^\alpha u(t) = a(t) u(t), \quad u(0) = y_2(0) - y_1(0) \ge 0,  \, u'(0) = y_2'(0) - y_1'(0)
	\end{equation}
	on $[0, T^*]$, and Theorem \ref{thm:no-intersect} or \ref{thm:no-intersect-0} 
	asserts that $u(t) > 0$ for all $t \in (0, T^*]$.
	
	Now we pick an arbitrary value $\tilde T \in (0,T^*]$. Since $[0, \tilde T] \subseteq [0, T^*]$,
	we can see that \eqref{eq:ivp-u} is satisfied on $[0, \tilde T]$. For $t$ in this interval, we then rewrite 
	the problem \eqref{eq:ivp-u} as
	\begin{align*}
		D^\alpha u(t) &= a^*(\tilde T) u(t) + (a(t) - a^*(\tilde T)) u(t), \\
		u(0) & = y_2(0) - y_1(0),  \qquad u'(0) = y_2'(0) - y_1'(0).
	\end{align*}
	By Theorem \ref{thm:voc}, it follows that
	\begin{align}
		\label{eq:voc-u1}
		u(t) & = (y_2(0) - y_1(0)) E_{\alpha,1}(a^*(\tilde T) t^\alpha)
				+ (y_2'(0) - y_1'(0)) t E_{\alpha,2}(a^*(\tilde T) t^\alpha) \\
					\nonumber
			& \qquad {} + \int_0^t (t - \tau)^{\alpha-1} E_{\alpha, \alpha}(a^*(\tilde T) (t - \tau)^\alpha)
					(a(\tau) - a^*(\tilde T)) u(\tau) \, \mathrm d \tau.
	\end{align}
	The definition of the function $a^*$ implies that $a(\tau) - a^*(\tilde T) \le 0$ for all 
	$\tau \in [0, T^*] \supseteq [0, t]$. In combination with the positivity of $u$ in the relevant
	interval that we had observed at the beginning of the proof and the choice of $T^*$
	which implies that the Mittag-Leffler function inside the integral 
	on the right-hand side of \eqref{eq:voc-u1} has a positive function value, 
	we conclude that this integral is nonnegative for all $t \in [0, \tilde T]$ and hence,
	in particular, for $t = \tilde T$. Therefore, we obtain
	\[
		|u(\tilde T)| = u(\tilde T) 
			\le (y_2(0) - y_1(0)) E_{\alpha,1}(a^*(\tilde T) \tilde T^\alpha)
				+ (y_2'(0) - y_1'(0)) \tilde T E_{\alpha,2}(a^*(\tilde T) \tilde T^\alpha) .
	\]
	Since $\tilde T$ was an arbitrary element of $(0, T^*]$,
	the upper bound for $|y_2(t) - y_1(t)|$ claimed in eq.~\eqref{eq:sep-linear} follows
	from the given assumption \eqref{eq:tstar1}
	upon the change of variable $\tilde T \mapsto t$.
	
	To prove the lower bound, we also start with the problem \eqref{eq:ivp-u}
	and again pick some arbitrary $\tilde T \in (0,T^*]$. Then we consider the initial 
	value problem
	\begin{align*}
		D^\alpha v(t) &= a_*(\tilde T) v(t), \\
		v(0) & = y_2(0) - y_1(0),  \qquad v'(0) = y_2'(0) - y_1'(0),
	\end{align*}
	for $t \in [0, \tilde T]$.
	According to \cite[Theorem 7.2]{Diethelmbook}, the unique solution to this problem is
	\begin{equation}
		\label{eq:v}
		v(t) = (y_2(0) - y_1(0)) E_{\alpha,1}(a_*(\tilde T) t^\alpha)
				+ (y_2'(0) - y_1'(0)) t E_{\alpha,2}(a_*(\tilde T) t^\alpha) .
	\end{equation}
	Then, for $t \in [0, \tilde T]$, we set $h(t) := u(t) - v(t)$ and see that the function $h$
	is the unique solution of the initial value problem
	\[
		D^\alpha h(t) = a_*(\tilde T) h(t) + (a(t) - a_*(\tilde T)) u(t), \qquad
		h(0) = h'(0) = 0,
	\]
	and so Theorem \ref{thm:voc} implies that
	\[
		h(t) = \int_0^t (t - \tau)^{\alpha - 1} E_{\alpha, \alpha}(a_*(\tilde T) (t-\tau)^\alpha)
				(a(\tau) - a_*(\tilde T)) u(\tau) \, \mathrm d \tau.
	\]
	From this relation, we can see (using the same arguments as above) that $h(t) \ge 0$
	for all $t \in [0, T^*]$, and so $u(t) \ge v(t)$ for all $t$. Hence, by definition of $u(t)$
	and \eqref{eq:v}, the lower bound for $|y_2(t) - y_1(t)| $ stated in \eqref{eq:sep-linear} 
	follows as well. 
	\qed
\end{pf}

\begin{rmk}
	Much as in \cite[Remark 2]{DT:separation1}, we can see that the upper bound in \eqref{eq:sep-linear}
	also holds if the differential equation under consideration
	is vector valued, i.e.\ if $f : [0, T] \times \mathbb R^d \to \mathbb R^d$ and 
	$y : [0,T] \to \mathbb R^d$ with arbitrary $d \in \mathbb N$, whereas the
	lower bound holds only in the scalar case ($d = 1$).
\end{rmk}

\begin{rmk}
	The observations of Remark \ref{rmk:cond-signs} apply to Theorem \ref{thm:sep-linear} as well.
\end{rmk}

\subsection{Main results for the nonlinear case}
\label{subs:nonlinear}

Based on Theorem \ref{thm:sep-linear}, we can now establish similar results for the 
general nonlinear case. In this situation, we temporarily impose a restriction on the given 
initial value problem \eqref{eq:ivp}, namely we assume that the function 
$f$ on the right-hand side of the differential equation has the property $f(t, 0) = 0$
for all $t \in [0,T]$. Our first results in this subsection will require this limitation; 
later on, we will remove it.

In complete analogy with the notions introduced in Theorem \ref{thm:sep-linear},
we start our considerations in the nonlinear case (under the restriction mentioned above) by defining
\begin{equation}
	\label{eq:astar}
	a_*(t) := \inf_{\tau \in [0, t], x \in \mathbb R \setminus \{0\}} \frac{f(\tau, x)} x
	\quad \mbox{ and } \quad
	a^*(t) := \sup_{\tau \in [0, t], x \in \mathbb R \setminus \{0\}} \frac{f(\tau, x)} x .
\end{equation}
Clearly, in the linear case $f(t, y) = a(t) y$ discussed in Subsection \ref{subs:linear}, these
definitions coincide with those introduced in the formulation of Theorem \ref{thm:sep-linear}.
Moreover, in the same way as in \cite[Proof of Lemma 2]{DT:separation1}, we can see
that these quantities are well defined and that the functions $a_*$ and $a^*$ are continuous
on $[0,T]$ if the function $f$ satisfies a Lipschitz condition.

Our first result then exploits the fact that, under the assumption $f(t, 0) = 0$ for all $t$,
the function given by $y_1(t) = 0$ for all $t$ is a solution to the differential equation given 
in \eqref{eq:ivp}.

\begin{thm}
	\label{thm:sep-nonlinear}
	Consider the differential equation $D^\alpha y_2(t) = f(t, y_2(t))$ for $t \in [0,T]$
	with $1 < \alpha < 2$ where $f$ is assumed to be continuous and to satisfy a Lipschitz condition in the
	second variable with Lipschitz constant $L$ and where $f(t, 0) = 0$ for all $t$.
	Additionally, let $T^* \in (0, T]$ satisfy the conditions \eqref{eq:tstar2}, \eqref{eq:tstar1a2} 
	and \eqref{eq:tstar1a3}.
	Then, the solution $y_2$ to this differential equation subject to the initial conditions
	$y_2(0) = y_{2,0}$ and $y_2'(0) = y_{2,1}$ satisfies the following estimates 
	for all $t \in [0,T^*]$:
	\begin{enumerate}[(a)]
	\item If $y_{2,0} \ge 0$ and $y_{2,1} \ge 0$ then
		\begin{align*}
			\MoveEqLeft{y_2(0) E_{\alpha,1}(a_*(t) t^\alpha) + y_2'(0) t E_{\alpha,2}(a_*(t) t^\alpha)} \\
			& \le  y_2(t) 
				\le y_2(0) E_{\alpha,1}(a^*(t) t^\alpha) + y_2'(0) t E_{\alpha,2}(a^*(t) t^\alpha).
		\end{align*}
	\item If $y_{2,0} \le 0$ and $y_{2,1} \le 0$ then
		\begin{align*}
			\MoveEqLeft{y_2(0) E_{\alpha,1}(a^*(t) t^\alpha) + y_2'(0) t E_{\alpha,2}(a^*(t) t^\alpha)} \\
			& \le  y_2(t) 
				\le y_2(0) E_{\alpha,1}(a_*(t) t^\alpha) + y_2'(0) t E_{\alpha,2}(a_*(t) t^\alpha).
		\end{align*}
	\end{enumerate}
	Moreover, if $y_{2,0} = 0$ then these statements also hold if \eqref{eq:tstar1a2} is not
	satisfied.
\end{thm}

\begin{pf}
	If $y_{2,0} = y_{2,1} = 0$ then all statements are trivial, so we only need 
	to explicitly consider the case that at least one of these values is nonzero.
	Furthermore, the statements are also trivial for $t = 0$. Thus we only have to 
	show the inequalities for $t \in (0, T^*]$.
	
	For the proof of (a), we note that, as indicated above, our assumptions 
	imply that the function $y_1(t) = 0$ is the unique solution to the given differential equation 
	subject to the initial conditions $y_1(0) = y_1'(0) = 0$. 
	Thus, $y_2(t) > y_1(t) = 0$ for all $t \in (0, T^*]$ by Theorem \ref{thm:no-intersect} if $y_{2,0} > 0$ 
	and by Theorem \ref{thm:no-intersect-0} if $y_{2,0} = 0$. 
	
	Now, as in the proof of Theorem \ref{thm:sep-linear},
	choose an arbitrary $\tilde T \in (0, T^*]$. Then, we can see that, for $t \in [0, \tilde T]$,
	the function $y_2$ satisfies the differential equation
	\[
		D^\alpha y_2(t) = a_*(\tilde T) y_2(t) + g_*(t, y_2(t))
	\]
	with
	\[
		g_*(t, y_2) := \left[ \frac{f(t, y_2)} {y_2} - a_*(\tilde T) \right] y_2,
	\]
	and by definition of $a_*(\tilde T)$, the expression in square brackets on the right-hand side
	of the definition of $g_*(t, y_2)$ is nonnegative. Thus, by Theorem \ref{thm:voc},
	\begin{align*}
		y_2(t) & = y_{2,0} E_{\alpha,1}(a_*(\tilde T) t^\alpha)
				+ y_{2,1} t E_{\alpha,2}(a_*(\tilde T) t^\alpha) \\
					\nonumber
			& \qquad {} + \int_0^t (t - \tau)^{\alpha-1} E_{\alpha, \alpha}(a_*(\tilde T) (t - \tau)^\alpha)
					g_*(\tau, y_2(\tau)) \, \mathrm d \tau,
	\end{align*}
	and in view of our considerations above---in particular also taking into account that $y_2(t) > 0$ 
	for all $t \in (0, T^*]$---we see that the integral in this equation is nonnegative. Thus, for $t = \tilde T$ we find
	\[
		y_2(\tilde T) \ge y_{2,0} E_{\alpha,1}(a_*(\tilde T) \tilde T^\alpha) 
					+ y_{2,1} \tilde T E_{\alpha,2}(a_*(\tilde T) \tilde T^\alpha) 
	\]
	for arbitrary $\tilde T \in (0,T^*]$ (and also, trivially, for $\tilde T = 0$). This is equivalent to
	the lower bound in the claim of part (a).
	
	To show the upper bound in part (a), we note that $y_2$ also satisfies the differential equation
	\[
		D^\alpha y_2(t) = a^*(\tilde T) y_2(t) + g^*(t, y_2(t))
	\]
	with
	\[
		g^*(t, y_2) := \left[ \frac{f(t, y_2)} {y_2} - a^*(\tilde T) \right] y_2,
	\]
	and by definition of $a^*(\tilde T)$, the expression in square brackets on the right-hand side
	of the definition of $g^*(t, y_2)$ is nonpositive. Thus, here we have by Theorem~\ref{thm:voc}
	that
	\begin{align*}
		y_2(t) & = y_{2,0} E_{\alpha,1}(a^*(\tilde T) t^\alpha)
				+ y_{2,1} t E_{\alpha,2}(a^*(\tilde T) t^\alpha) \\
					\nonumber
			& \qquad {} + \int_0^t (t - \tau)^{\alpha-1} E_{\alpha, \alpha}(a^*(\tilde T) (t - \tau)^\alpha)
					g^*(\tau, y_2(\tau)) \, \mathrm d \tau
	\end{align*}
	where now the integral is nonpositive, so by an argument as above we also establish the upper bound.
	
	Statement (b) can be shown in a completely analog manner.
	\qed
\end{pf}

As announced above, we now finish our quantitative discussion of upper and lower bounds for the 
difference between two distinct solutions to the same fractional differential equation of order $\alpha \in (1,2)$
by a look at the general nonlinear case without assuming the hypothesis that $f(t,0) = 0$ for all $t$.
Note that this final theorem of this section is the only theorem that requires knowledge of one
specific solution to the differential equation in question. All other theorems proved so far in this paper
were exclusively formulated in terms of the given function $f$ on the right-hand side of the differential
equation and the given initial values.

\begin{thm}
	Consider the initial value problems \eqref{eq:ivp} for $t \in [0,T]$
	with $1\! < \! \alpha \! <\! 2$ where $f$ is assumed to be continuous and to satisfy a Lipschitz condition in the
	second variable with Lipschitz constant $L$. Assume that $y_{1,0} \le y_{2,0}$ and $y_{1,1} \le y_{2,1}$.
	Moreover, let $T^* \in (0, T]$ satisfy the conditions \eqref{eq:tstar2}, \eqref{eq:tstar1a2} and \eqref{eq:tstar1a3}.
	Then, defining
	\begin{align*}
		\tilde a_*(t) & := \inf_{\tau \in [0,t], y \ne 0} \frac{f(\tau, y + y_1(\tau)) - f(\tau, y_1(\tau))} y
	\intertext{and}
		\tilde a^*(t) & := \sup_{\tau \in [0,t], y \ne 0} \frac{f(\tau, y + y_1(\tau)) - f(\tau, y_1(\tau))} y,
	\end{align*}
	we have
	\begin{align*}
		\MoveEqLeft{(y_{2,0} - y_{1,0}) E_{\alpha,1}(\tilde a_*(t) t^\alpha) 
					+ (y_{2,1} - y_{1,1}) t E_{\alpha,2}(\tilde a_*(t) t^\alpha)} \\
		& \le  y_2(t) - y_1(t) \\
		& \le (y_{2,0} - y_{1,0}) E_{\alpha,1}(\tilde a^*(t) t^\alpha) 
					+ (y_{2,1} - y_{1,1}) t E_{\alpha,2}(\tilde a^*(t) t^\alpha)
	\end{align*}
	for all $t \in [0, T^*]$.
\end{thm}

\begin{pf}
	Defining
	\[
		\tilde f(t, y) := f(t, y + y_1(t)) - f(t, y_1(t)),
	\]
	we see the following facts:
	\begin{enumerate}[(a)]
	\item The function $u := y_2 - y_1$ satisfies the differential equation
		\begin{equation}
			\label{eq:ivp-nl}
			D^\alpha u (t) 
			= D^\alpha y_2(t) - D^\alpha y_1(t)
			= f(t, y_2(t)) - f(t, y_1(t))
			= \tilde f(t, u(t))
		\end{equation}
		and the initial conditions $u(0) = y_{2,0} - y_{1,0}$ and $u'(0) = y_{2,1} - y_{1,1}$.
	\item $\tilde f(t, 0) = 0$ for all $t$.
	\item $\tilde f$ satisfies a Lipschitz condition with the same Lipschitz constant $L$ as $f$.
	\end{enumerate}
	Therefore, we may apply Theorem \ref{thm:sep-nonlinear}(a) to the solution $u$ of the differential equation 
	\eqref{eq:ivp-nl} and derive the desired result.
	\qed
\end{pf}

\section{Zeros of Mittag-Leffler functions}
\label{sec:ml-zeros}

The conditions \eqref{eq:tstar2}, \eqref{eq:tstar1a2} and \eqref{eq:tstar1a3} play a significant role in deciding when
the results of Section \ref{sec:separation} can be applied. To obtain a thorough understanding
of the meaning of these conditions, it is essential to have some knowledge of the locations
of the smallest positive zeros of the three functions 
\[
	z \mapsto E_{\alpha, \beta}(-z^\alpha) 
	\qquad \text{ with } \beta \in \{1, \alpha, 2\}.
\]
While a large amount of information is available on the zeros with large
modulus of these functions (see, e.g., \cite[Section 4.6]{GKMR2020}), we have unfortunately 
been unable to find any 
concrete statements about the zeros with small modulus. Therefore, to allow the reader to at least obtain
a reasonable impression of the situation, we now present some numerical results. To this end, 
we define 
\begin{equation}
	\label{eq:minzero}
	Z_{\beta}(\alpha) := \min \{ z > 0 : E_{\alpha, \beta}(-z^\alpha)  = 0 \}
	\quad \text{ for } \beta \in \{1, \alpha, 2\}.
\end{equation}
Throughout this section, whenever we had to evaluate some Mittag-Leffler function, we have used
Garrappa's MATLAB implementation of his algorithm described in \cite{Ga2015b}.

We hope that our paper will initiate further analytical investigations on the question for the exact
values of $Z_\beta(\alpha)$.

\subsection{Small positive zeros of $z \mapsto E_{\alpha, \alpha}(-z^\alpha)$}
\label{sec:ml-zeros-alpha}

\begin{figure}[htp]
	\centering
	\includegraphics[width=0.75\textwidth]{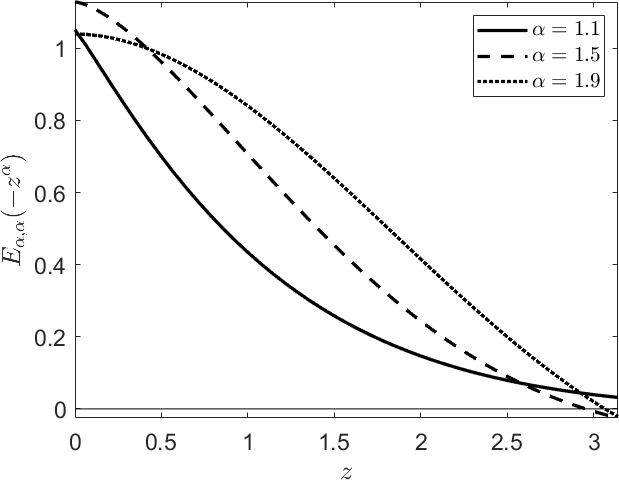}
	\caption{\label{fig:alpha-1} Plots of $E_{\alpha, \alpha}(-z^\alpha)$ vs.\ $z$ for $\alpha \in \{ 1.1, 1.5, 1.9 \}$.}
\end{figure}

We begin our discussion of the quantity $Z_\beta(\alpha)$ defined in \eqref{eq:minzero}
with the case $\beta = \alpha$. From \cite[Section 7]{HAV2013} we obtain that the number of
zeros of this function is positive and finite.
A limit consideration using the power series representation of the Mittag-Leffler function yields 
\begin{subequations}
\begin{align}
	\label{eq:alpha1}
	\lim_{\alpha \to 1+} E_{\alpha, \alpha}(-z^\alpha) & = E_{1,1}(-z) = \exp(-z) \\
\intertext{and}
	\label{eq:alpha2}
	\lim_{\alpha \to 2-} E_{\alpha, \alpha}(-z^\alpha) & = E_{2,2}(-z^2) = \frac 1 z \sin z.
\end{align}
\end{subequations}
From \eqref{eq:alpha2} we deduce that the smallest positive $z$ for which $E_{2,2}(-z^2) = 0$
is $z = \pi$, and so we have $Z_2(2) = \pi$. On the other hand, by \eqref{eq:alpha1}, the equation 
$E_{1,1}(-z) = 0$ has no positive solution at all, and thus we may expect 
$\lim_{\alpha \to 1+} Z_\alpha(\alpha) = \infty$. For $\alpha \in (1,2)$, it is a relatively easy 
matter to compute the values of $Z_\alpha(\alpha)$ numerically: The graphs of the function (see the examples
given in Figure \ref{fig:alpha-1}) indicate that we may use the standard
Newton iteration algorithm with starting value $z_0 = \pi/2$; indeed the numerical experiments reveal that
this always leads to a rapid convergence against the correct solution. To be precise, the number of iterations required
to compute the first positive zero of $E_{\alpha, \alpha}(-z^\alpha)$ with an absolute error of $10^{-12}$
was 19 for $\alpha = 1.001$, decreased monotonically to 10 as $\alpha$ increased to $1.045$, further to 7
as $\alpha$ reached the value $1.303$, and stayed at or below 6 for $\alpha \ge 1.615$.
The required derivative of the function 
$E_{\alpha, \alpha}(-z^\alpha)$ can be computed using the chain rule and a classical formula
for derivatives of Mittag-Leffler functions \cite[eq.~(4.3.2)]{GKMR2020}, thus obtaining
\[
	\frac{\mathrm d}{\mathrm d z} E_{\alpha, \alpha}(-z^\alpha) 
	= \frac 1 z \left( E_{\alpha, \alpha-1}(-z^\alpha) + (1-\alpha) E_{\alpha, \alpha}(-z^\alpha) \right).
\]

\begin{figure}[htp]
	\includegraphics[width=0.49\textwidth]{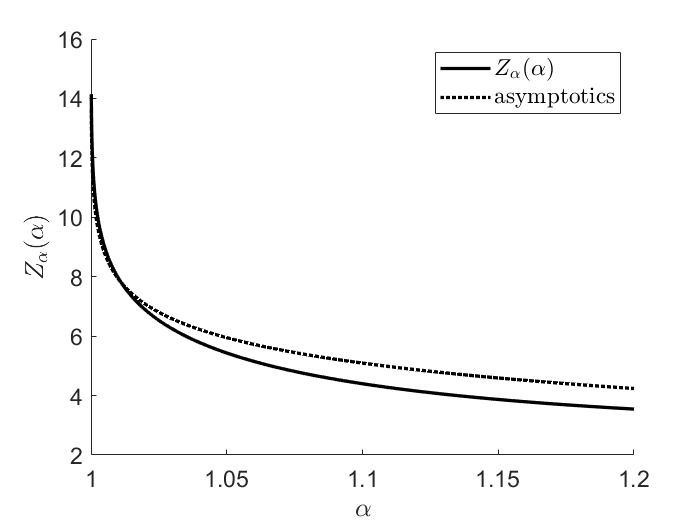}
	\hfill
	\includegraphics[width=0.49\textwidth]{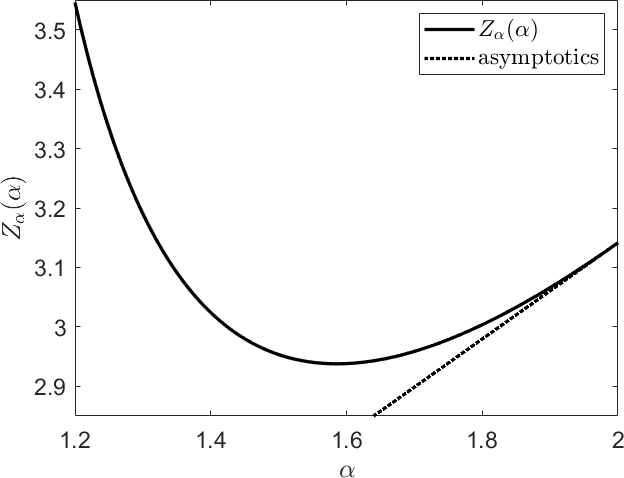}
	\caption{\label{fig:alpha-2} Left: Location $Z_\alpha(\alpha)$ of smallest positive zero of 
		$E_{\alpha, \alpha}(-z^\alpha)$ vs.\ $\alpha$ for $\alpha \in (1, 1.2]$
		(continuous line) and suspected asymptotic behaviour for $\alpha \to 1+$ (dotted line); 
		Right: Location $Z_\alpha(\alpha)$ of smallest positive zero of $E_{\alpha, \alpha}(-z^\alpha)$ vs.\ $\alpha$ 
		for $\alpha \in (1.2, 2]$
		(continuous line) and suspected asymptotic behaviour for $\alpha \to 2-$ (dotted line).}
\end{figure}

The algorithm then yields the numerical results shown in the continuous line in the graphs of Figure \ref{fig:alpha-2}.
To provide a clearer picture, each of the two graphs in this figure covers only a part of the interval $\alpha \in (1,2]$.

Specifically, we find from our numerical results that $\min_{\alpha \in (1,2]} Z_\alpha(\alpha) \approx 2.9378538$ and
that this minimum is attained for $\alpha \approx 1.586$. In particular, $Z_\alpha(\alpha)$ does not depend on $\alpha$
in a monotonic way.

Also, we can analyze the asymptotic behaviour of 
$Z_{\alpha}(\alpha)$ for $\alpha \to 1+$ and $\alpha \to 2-$, respectively, based on the numerical results. 
The observations obtained in this way give rise to the following conjecture. To support items (b) and (c) of the 
conjecture, the suspected asymptotics are plotted as the dotted lines in the two parts of Figure \ref{fig:alpha-2}.

\begin{conj}
	\begin{enumerate}[(a)]
	\item For any $\alpha \in (1,2]$, the statement $E_{\alpha, \alpha}(-z^\alpha) > 0$ holds for all $z \in [0, 2.93785]$.
	\item For $\alpha \to 1+$, we have $Z_\alpha(\alpha) = (c \ln(\alpha-1) + d) (1 + o(1))$ 
		with $c \approx -0.81$ and $d \approx 2.25$.
	\item For $\alpha \to 2-$, we have $Z_\alpha(\alpha) = \pi + c (2 - \alpha)  + o(2 - \alpha)$ with $c \approx -0.81$.
	\end{enumerate}
\end{conj}

Note that coefficients with values approximately equal to $-0.81$ occur both in the asymptotic
expansion for $\alpha \to 1+$ and in the expansion for $\alpha \to 2-$. We do not know whether
or not these two coefficients are exactly identical.

\subsection{Small positive zeros of $z \mapsto E_{\alpha, 2}(-z^\alpha)$}
\label{sec:ml-zeros-2}

\begin{figure}[htp]
	\centering
	\includegraphics[width=0.75\textwidth]{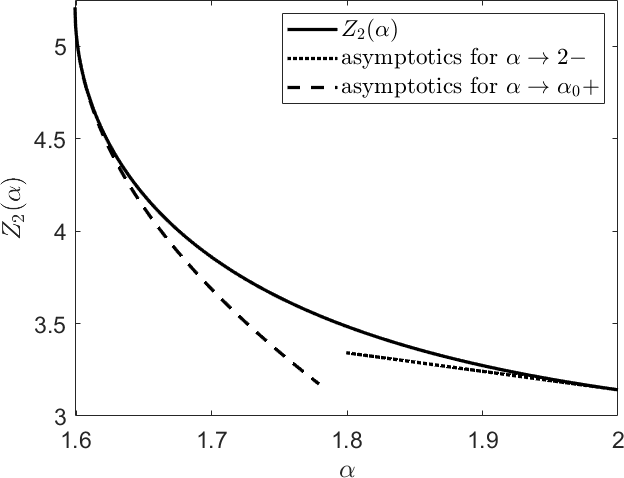}
	\caption{\label{fig:2-1} Location $Z_2(\alpha)$ of smallest positive zero of 
		$E_{\alpha, 2}(-z^\alpha)$ vs.\ $\alpha$ for $\alpha \in [1.59911520635, 2]$
		(continuous line) and suspected asymptotic behaviour according to Conjecture \ref{conj:2} 
		for $\alpha \to 2-$ (dotted line) and for $\alpha \to \alpha_0+$ (dashed line).}
\end{figure}

Next, we look at the zeros of $E_{\alpha, \beta}(-z^\alpha)$ for $\beta = 2$. Here, we observe
a substantially different behaviour: The function has real zeros if and only if $\alpha$ is above or equal to a certain 
threshold $\alpha_0$ whose value is not known exactly \cite[Sections 6 and 7]{HAV2013}.
The approximation provided in \cite[eq.~(8)]{HAV2013} indicates that $\alpha_0 \approx 1.59927$; 
our numerical observations indicate that the correct value is $\alpha_0 \approx 1.599115206$ (see the right part of Figure
\ref{fig:2-2}).
We can then actually proceed in the way indicated for the case $\beta = \alpha$
in Subsection \ref{sec:ml-zeros-alpha} for $\alpha \in [\alpha_0, 2]$ while
for $\alpha \in (1, \alpha_0)$ the search for zeros is pointless.
We can summarize our findings as follows (see, in particular, Figures \ref{fig:2-2} and \ref{fig:2-1} 
for a motivation of these formulations):

\begin{figure}[htp]
	\includegraphics[width=0.49\textwidth]{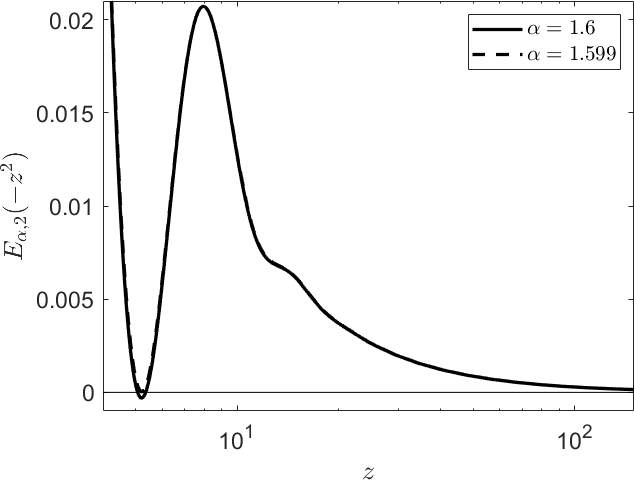}
	\hfill
	\includegraphics[width=0.49\textwidth]{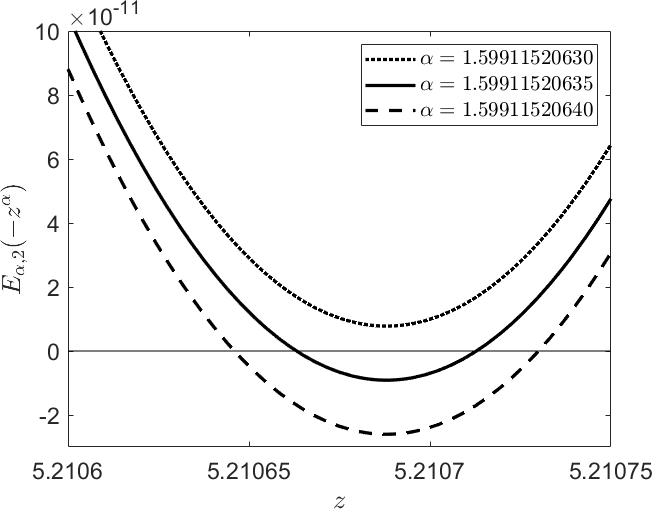}
	\caption{\label{fig:2-2} Plots of $E_{\alpha, 2}(-z^\alpha)$ vs.\ $z$ for $\alpha = 1.599$ 
		(dashed line) and $\alpha = 1.6$ (continuous line).
		Left: Plot range $z \in [4, 150]$ with logarithmically scaled horizontal axis showing that both functions 
		(visually hard to distinguish from each other) 
		seem to remain positive for large~$z$; 
		Right: Zoom into $z \in [5.2106, 5.21075]$ with linear scale on horizontal axis 
		reveals that the function for $\alpha = 1.59911520635$ has a zero at $z \approx 5.210663$
		whereas the function for $\alpha = 1.5991152064$ has a local minimum with a strictly positive function value 
		and hence does not have a zero.}
\end{figure}

\begin{conj}
	\label{conj:2}
	There exists some $\alpha_0 \in (1.59911, 1.59912)$ such that the following statements hold:
	\begin{enumerate}[(a)]
	\item For each $\alpha \in [\alpha_0, 2]$, the equation $E_{\alpha, 2}(-z^\alpha) = 0$ possesses at least one solution in $(0, \infty)$.
		Denoting (as in eq.~\eqref{eq:minzero} above) the smallest of these solutions by $Z_2(\alpha)$, we have:
		\begin{enumerate}[(i)]
		\item $Z_2(\alpha) = \pi + (2-\alpha) + o(2 - \alpha)$ for $\alpha \to 2-$, 
		\item $Z_2(\alpha) = Z_2(\alpha_0) + c (\alpha - \alpha_0)^b (1 + o(1))$ 
			for $\alpha \to \alpha_0+$ where $Z_2(\alpha_0) \approx 5.21066$,
			$b \approx 1/2$ and $c \approx -4.8$,
		\item $Z_2(\alpha)$ is a strictly decreasing function of $\alpha$.
		\end{enumerate}
	\item For each $\alpha \in (1, \alpha_0)$, we have $E_{\alpha, 2}(-z^\alpha) > 0$ for all $z \ge 0$.
	\end{enumerate}
\end{conj}

\subsection{Small positive zeros of $z \mapsto E_{\alpha, 1}(-z^\alpha)$}
\label{sec:ml-zeros-1}

\begin{figure}[htp]
	\centering
	\includegraphics[width=0.75\textwidth]{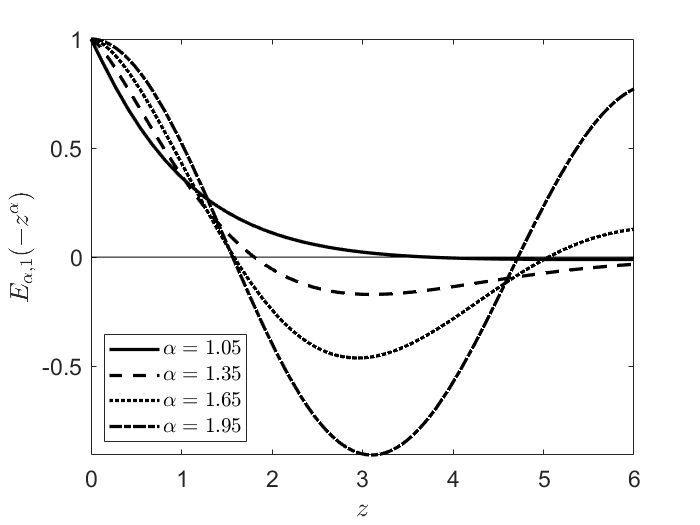}
	\caption{\label{fig:1-1} Plots of $E_{\alpha, 1}(-z^\alpha)$ vs.\ $z$ for $\alpha \in \{ 1.05, 1.35, 1.65, 1.95 \}$.}
\end{figure}

As far as the function $E_{\alpha, 1}(-z^\alpha)$ is concerned, 
we once again obtain from \cite{HAV2013} that it has a positive and finite
number of zeros for any $\alpha \in (1,2)$.
Since
\begin{subequations}
\begin{align}
	\label{eq:11}
	\lim_{\alpha \to 1+} E_{\alpha, 1}(-z^\alpha) & = E_{1,1}(-z) = \exp(-z) \\
\intertext{and}
	\label{eq:12}
	\lim_{\alpha \to 2-} E_{\alpha, 1}(-z^\alpha) & = E_{2,1}(-z^2) = \cos z,
\end{align}
\end{subequations}
arguments similar to those used above yield $Z_1(2) = \pi / 2$. On the other hand, by \eqref{eq:11}, the equation 
$E_{1,1}(-z) = 0$ has no positive solution at all, and thus we may expect 
$\lim_{\alpha \to 1+} Z_1(\alpha) = \infty$. The graphs of these functions for some values of $\alpha$ that cover
a wide range of our interval of interest are shown in Figure~\ref{fig:1-1}. The shapes of these graphs indicate 
that we may find the smallest zeros of each of these functions with the help of a Newton iteration with starting
value $\pi/2$, and our numerical experiments confirm this expectation. Indeed, for all $\alpha \in [1.001, 2]$ 
the method converged up to an accuracy of $10^{-12}$ in at most 13 iterations.

\begin{figure}[htp]
	\includegraphics[width=0.49\textwidth]{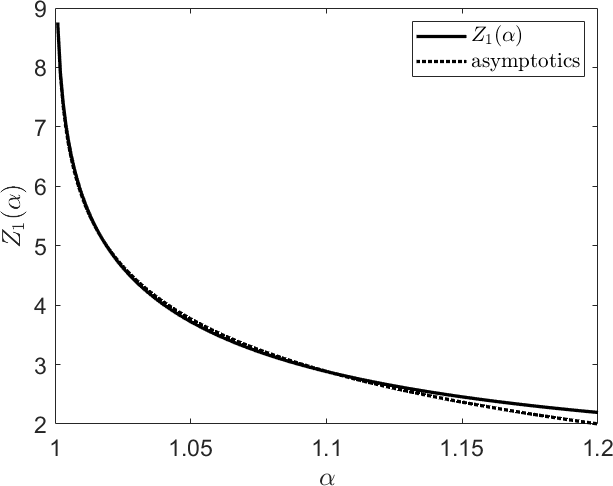}
	\hfill
	\includegraphics[width=0.49\textwidth]{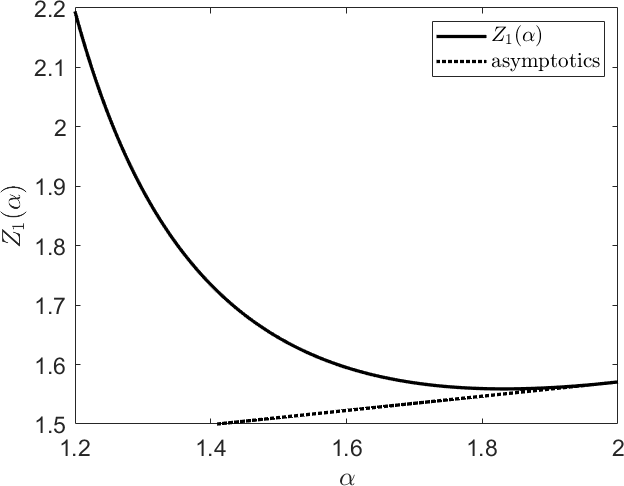}
	\caption{\label{fig:1-2} Left: Location $Z_1(\alpha)$ of smallest positive zero of 
		$E_{\alpha, 1}(-z^\alpha)$ vs.\ $\alpha$ for $\alpha \in (1, 1.2]$
		(continuous line) and suspected asymptotic behaviour for $\alpha \to 1+$ (dotted line); 
		Right: Location $Z_1(\alpha)$ of smallest positive zero of $E_{\alpha, 1}(-z^\alpha)$ vs.\ $\alpha$ 
		for $\alpha \in (1.2, 2]$
		(continuous line) and suspected asymptotic behaviour for $\alpha \to 2-$ (dotted line).}
\end{figure}

The graphs in Figure \ref{fig:1-2} show the locations $Z_1(\alpha)$ of the smallest positive zero 
of $E_{\alpha, 1}(-z^\alpha)$ vs.\ $\alpha$ obtained in this way. In particular, we once agan see
(like in the case $\beta = \alpha$) that $Z_1(\alpha)$ does not depend on $\alpha$ monotonically.
More precisely, the smallest value of $Z_1(\alpha)$ is attained for $\alpha \approx 1.833$, with
$\min_{\alpha \in (1,2]} Z_1(\alpha) \approx 1.559$.

More precisely, our numerical calculations give rise to the following conjecture (see also Figures \ref{fig:1-1}
and \ref{fig:1-2}):

\begin{conj}
	\begin{enumerate}[(a)]
	\item For any $\alpha \in (1,2]$, the statement $E_{\alpha, 1}(-z^\alpha) > 0$ holds for all $z \in [0, 1.559]$.
	\item For $\alpha \to 1+$, we have $Z_1(\alpha) = (c \ln(\alpha-1) + d) (1 + o(1))$ 
		with $c \approx -0.275$ and $d \approx 17.54$.
	\item For $\alpha \to 2-$, we have $Z_1(\alpha) = \frac \pi 2 + c (2 - \alpha)  + o(2 - \alpha)$ with $c \approx -0.12$.
	\end{enumerate}
\end{conj}

\section{Conclusion and outlook}
\label{sec:conclusion}

As indicated in Section \ref{sec:intro}, we believe that our findings can be useful to
obtain an improved understanding of fractional boundary value problems of the form \eqref{eq:bvp}
and in the development and analysis of associated numerical solution techniques, in particular 
for algorithms based on the principle of shooting methods \cite{Di:shooting-bvp,Keller}, thus extending the
developments for terminal value problems of order $\alpha \in (0,1)$ described in \cite{DU:shooting1}.

To fully exploit the results presented here, it will be necessary to obtain further, and in particular
analytic, information 
about the location of the negative real zeros of the Mittag-Leffler functions used above. We hope that
our work will motivate the initiation of such investigations.

\section*{Acknowledgment}

The authors would like to thank Roberto Garrappa for the useful discussion about the questions
addressed in Section \ref{sec:ml-zeros} and, in particular, for pointing out reference \cite{HAV2013}.


\begin{thebibliography}{10}
\expandafter\ifx\csname url\endcsname\relax
  \def\url#1{\texttt{#1}}\fi
\expandafter\ifx\csname urlprefix\endcsname\relax\def\urlprefix{URL }\fi
\expandafter\ifx\csname href\endcsname\relax
  \def\href#1#2{#2} \def\path#1{#1}\fi

\bibitem{DT:separation1}
K.~Diethelm, H.~T. Tuan, Upper and lower estimates for the separation of
  solutions to fractional differential equations, Fract. Calc. Appl. Anal. 25
  (2022) 166--180.
\newblock \href {https://doi.org/10.1007/s13540-021-00007-x}
  {\path{doi:10.1007/s13540-021-00007-x}}.

\bibitem{DF}
K.~Diethelm, N.~J. Ford, Analysis of fractional differential equations, J.
  Math. Anal. Appl. 265 (2002) 229--248.
\newblock \href {https://doi.org/10.1006/jmaa.2000.7194}
  {\path{doi:10.1006/jmaa.2000.7194}}.

\bibitem{Diethelmbook}
K.~Diethelm, The Analysis of Fractional Differential Equations, Vol. 2004 of
  Lecture Notes in Mathematics, Springer, Berlin, 2010.
\newblock \href {https://doi.org/10.1007/978-3-642-14574-2}
  {\path{doi:10.1007/978-3-642-14574-2}}.

\bibitem{GKMR2020}
R.~Gorenflo, A.~Kilbas, F.~Mainardi, S.~Rogosin, Mittag-Leffler Functions,
  Related Topics and Applications, 2nd Edition, Springer, Berlin, 2020.
\newblock \href {https://doi.org/10.1007/978-3-662-61550-8}
  {\path{doi:10.1007/978-3-662-61550-8}}.

\bibitem{Tisdell2012}
C.~C. Tisdell, On the application of sequential and fixed-point methods to
  fractional differential equations of arbitrary order, J. Integral Equations
  Appl. 24 (2012) 283--319.
\newblock \href {https://doi.org/10.1216/JIE-2012-24-2-283}
  {\path{doi:10.1216/JIE-2012-24-2-283}}.

\bibitem{DU:shooting1}
K.~Diethelm, F.~Uhlig, A new approach to shooting methods for terminal value
  problems of fractional differential equations, J. Sci. Comput. 97 (2023).
\newblock \href {https://doi.org/10.1007/s10915-023-02361-9}
  {\path{doi:10.1007/s10915-023-02361-9}}.

\bibitem{Di:shooting-bvp}
K.~Diethelm, Shooting methods for fractional {Dirichlet}-type boundary value
  problems of order $\alpha \in (1,2)$ with {Caputo} derivatives, Tech. rep.,
  arXiv:2402.03487 (2024).

\bibitem{CT2017}
N.~D. Cong, H.~T. Tuan, Generation of nonlocal fractional dynamical systems by
  fractional differential equations, J. Integral Equations Appl. 29 (2017)
  585--608.
\newblock \href {https://doi.org/10.1216/JIE-2017-29-4-585}
  {\path{doi:10.1216/JIE-2017-29-4-585}}.

\bibitem{LV2008}
V.~Lakshmikantham, A.~S. Vatsala, Basic theory of fractional differential
  equations, Nonlinear Anal., Theory Methods Appl., Ser. A, Theory Methods 69
  (2008) 2677--2682.
\newblock \href {https://doi.org/10.1016/j.na.2007.08.042}
  {\path{doi:10.1016/j.na.2007.08.042}}.

\bibitem{Ga2015b}
R.~Garrappa, Numerical evaluation of two and three parameter {Mittag-Leffler}
  functions, SIAM J. Numer. Anal. 53 (2015) 1350--1369.
\newblock \href {https://doi.org/10.1137/140971191}
  {\path{doi:10.1137/140971191}}.

\bibitem{HAV2013}
J.~W. Hanneken, B.~N.~N. Achar, D.~M. Vaught, An alpha-beta phase diagram
  representation of the zeros and properties of the {M}ittag-{L}effler
  function, Adv. Math. Phys.\ 421685 (2013).
\newblock \href {https://doi.org/10.1155/2013/421865}
  {\path{doi:10.1155/2013/421865}}.

\bibitem{Keller}
H.~B. Keller, Numerical Methods for Two-Point Boundary-Value Problems, Dover,
  Mineola, 2018.

\end{thebibliography}
\end{document}